\documentclass[11pt]{article} 
\usepackage{amsmath, amssymb, amsthm,fancyhdr}

\newcommand{\version}{version 3.0,\ \   July 31, 2009}

\newcommand{\6}{\partial}
\newcommand{\arrow}{{\:\longrightarrow\:}}
\newcommand{\R}{{\Bbb R}}

\numberwithin{equation}{section}
\newtheorem{theorem}[equation]{Theorem}
\newtheorem*{theorem*}{Theorem}

\newtheorem{lemma}[equation]{Lemma}

\theoremstyle{remark}
\newtheorem{remark}[equation]{\bf Remark}
\newtheorem*{remark*}{\bf Remark}
\newtheorem{remarks}[equation]{\bf Remarks}
\newtheorem{definition}[equation]{\bf Definition}

\newtheorem{note*}{\bf Note}

\newtheorem{example*}{\bf Example}

\begin{document}

\begin{center}
{\LARGE\bf
Stable bundles on hypercomplex surfaces
}
\\[4mm]
Ruxandra Moraru, \footnote{R. Moraru is supported by an NSERC grant}
Misha Verbitsky,\footnote{M. Verbitsky is
supported by EPSRC grant  GR/R77773/01}\\[4mm]
{\tt moraru@math.uwaterloo.ca, verbit@mccme.ru}
\end{center}

{\small 
\hspace{0.15\linewidth}
\begin{minipage}[t]{0.8\linewidth}
{\bf Abstract} \\
A hypercomplex manifold is a manifold equipped with
three complex structures $I, J, K$  satisfying the quaternionic
relations. Let $M$ be a 4-dimensional compact smooth 
manifold equipped with a hypercomplex structure,
and $E$ be a vector bundle on $M$. We show that the 
moduli space of anti-self-dual connections on $E$ 
is also hypercomplex, and admits a strong HKT metric. 
We also study manifolds with (4,4)-supersymmetry,
that is, Riemannian manifolds  equipped with a pair 
of strong HKT-structures that have opposite torsion. 
In the language of Hitchin's and Gualtieri's 
generalized complex geometry, (4,4)-manifolds are 
called ``generalized hyperk\"ahler manifolds''.
We show that the moduli space of anti-self-dual 
connections on $M$  is a (4,4)-manifold if 
$M$ is equipped with a (4,4)-structure.
\end{minipage}
}

\tableofcontents


\section{Introduction}


\subsection{Instanton moduli and stable holomorphic bundles}

Ever since it was established by Donaldson and Uhlenbeck-Yau,
the correspondence between instantons and stable holomorphic vector bundles
on K\"ahler manifolds has been a constant source of new information 
about both instantons and holomorphic vector bundles. 

One of the most immediate applications of this correspondence is
the following. Let $(X,g)$ be a compact Riemannian 4-dimensional
manifold admitting complex structures $I_1$, $I_2$ such that the
metric $g$ is K\"ahler with respect to both $I_1$ and $I_2$.
This happens, for instance, when $X$ is a hyperk\"ahler
4-manifold, that is, a K3 surface or a compact
complex torus.

The moduli space $\cal M$ of instantons on $X$ depends only on a 
metric. From the Donaldson-Uhlenbeck-Yau
theorem, we obtain that $\cal M$, as a topological space,
is identified with the moduli of stable holomorphic bundles $E$
with $c_1(E)=0$ on $(X,I_1)$ and on $(X, I_2)$. Therefore,
$\cal M$ is equipped with a pair of complex
structures, one induced from $I_1$, another from $I_2$.

When $X$ is hyperk\"ahler, this can be used to show
that the moduli of instantons is hyperk\"ahler 
as well. This result was obtained by
A. Tyurin (\cite{_Tyurin:K3_}) and generalized in 
\cite{V-hyper} to hyperk\"ahler manifolds of
arbitrary dimension.

For non-K\"ahler complex manifolds, a version of the
Donald\-son\--Uhlen\-beck-\-Yau theorem was obtained by 
Buchdahl \cite{Buchdahl} for surfaces
and by Li and Yau \cite{LY} for general Hermitian
manifolds. In this context, the correspondence
between instantons and stable holomorphic
vector bundles is usually called {\em the 
Kobayashi-Hitchin correspondence}.
This result is not new, but the
full impact of the Buchdahl-Li-Yau theorem
in the geometry of non-K\"ahler manifolds is still 
not completely realized, though a book by L\"ubke and Teleman
(\cite{LT}) studies it in wonderful detail. 

Let $(X, I, g)$ be a compact complex Hermitian manifold,
and $\omega\in \Lambda^{1,1}(X)$ its Hermitian form.
If $\6\bar\6(\omega^{\dim X -1})=0$, then the Hermitian
metric on $X$ is called {\bf a Gauduchon metric}.
P. Gauduchon (\cite{_Gauduchon_1984_}) has proven 
that such a metric exists in each conformal class,
and is unique up to a constant. 

When $(X, I, g)$ is equipped with a Gauduchon metric,
the Li-Yau theorem identifies the instanton moduli
space with the space of stable holomorphic bundles
(see Section \ref{_Instanton_Section_} for details).

In this context, an instanton is a Hermitian bundle $E$
with a connection $A$ whose curvature 2-form $F_A$
is of type $(1,1)$ and is pointwise
orthogonal to the Hermitian form:
\[
F_A\in \Lambda^{1,1}(X), \ \ F_A \bot \omega.
\]
When $X$ is a complex surface, these conditions
are equivalent to the anti-self-duality of $A$
(see Section \ref{_Instanton_Section_}).

Now, assume that $(X,g)$ is a Riemannian manifold
admitting two complex structures $I_1$, $I_2$,
such that $g$ is Hermitian and Gauduchon with
respect to both $I_1$ and  $I_2$. Then the Buchdahl-Li-Yau
theorem implies that the moduli $\cal M$ 
of anti-seld-dual connections on $X$ is equipped 
with two complex structures, induced by $I_1$ and $I_2$. 
It is not generally known how these
complex structures relate to each other. However, if
$X$ is equipped with an additional geometric structure
(HKT- or bi-Hermitian), then it is possible
to recover a similar structure on the
 moduli space.

\subsection{Bismut connection and HKT-structures}

\begin{definition}
Let $(M,I,g)$ be a complex Hermitian manifold.
A connection $\nabla:\; TM\arrow TM \otimes \Lambda^1 M$
is called {\bf Hermitian} if $\nabla I = \nabla g=0$.
Consider its torsion $T_1\in \Lambda^2 \otimes TM$,
and let $T \in \Lambda^2 \otimes \Lambda^1M$
be the tensor obtained from $T_1$ via the
isomorphism $TM \cong \Lambda^1 M$ 
provided by $g$. A Hermitian connection
is called {\bf a Bismut connection}, or
{\bf a connection with skew-symmetric torsion},
if $T$ is {\it skew-symmetric}, that 
is, lies in 
$\Lambda^3 M \subset \Lambda^2 \otimes \Lambda^1M$.
The 3-form $T$ is called {\bf the torsion form}
of Bismut connection.
\end{definition}

\begin{theorem}\label{_Bismut_Theorem_}
Let $(M,I,g)$ be a complex Hermitian manifold.
Then $M$ admits a Bismut connection $\nabla$,
which is unique. Moreover, its torsion form
is equal to $Id\omega$.
\end{theorem}

\begin{proof} 
See
\cite{_Bismut:connection_},
\cite{_Friedrich_Ivanov:physics_}. 
\end{proof}

\begin{remark}
Clearly, if $d\omega=0$, then the Bismut connection
is torsion-free, and thus coincides with the Levi-Civita
connection. Theorem \ref{_Bismut_Theorem_}
can therefore be used to show that the Levi-Civita connection
on a K\"ahler manifold satisfies $\nabla I=0$.
\end{remark}

Connections with skew-symmetric torsion play an
important role in string physics (see for example
\cite{_Ivanov_Papadopoulos_}). In the physics
literature, a complex Hermitian 
manifold $(M,I,g)$ with a Bismut connection 
is called {\bf a KT-manifold} (K\"ahler torsion manifold).
If, in addition, the torsion 3-form is closed, then
$(M,I,g)$ is called {\bf a strong KT-manifold}.
By Theorem \ref{_Bismut_Theorem_},
a manifold is therefore strong KT if and only if
$\6\bar\6 \omega=0$. For complex surfaces, 
this is equivalent to $g$ being a Gauduchon 
metric.

There are several other structures based on Bismut
connections which are even more important.

\begin{definition}
Let $(M,g)$ be a Riemannian manifold, and $I_-$, $I_+$ be
Hermitian complex structures. Consider the corresponding
Bismut connections, and suppose that their torsion 3-forms
satisfy $T_+= - T_-$, $d T_\pm=0$.
Then $(M,g, I_+, I_-)$ is called {\bf bi-Hermitian}.
\end{definition}

Bi-Hermitian structures appear naturally in 
several different (and seemingly unrelated) contexts.
In differential geometry, these were studied
by Apostolov, Gauduchon and Grantcharov (\cite{_AGG_}),
who obtained classification results in the case when
$\dim_\R M=4$; they showed in particular that if
$M$ is of K\"ahler type,
then it is a rational surface, a torus or a K3 surface. 
In physics, such structures
were studied as early as 1984  by Gates, Hull and
Ro\v cek (\cite{_GHR_}),
in connection to $D=2$, $N=4$ supersymmetric
$\sigma$-models. More recently, bi-Hermitian manifolds
have appeared both in mathematics and in string physics,
due to the work of N. Hitchin and M. Gualtieri on
generalized complex geometry. In his Ph. D. thesis,
\cite{_Gualtieri:thesis_}, Gualtieri explored 
the notion of generalized complex manifold, which was
first developed by Hitchin (\cite{_Hitchin:gene_}). 
He defined generalized K\"ahler manifolds, and 
described them in terms of more classical
differential-geometric structures.
More precisely, Gualtieri
found that a generalized K\"ahler structure
on a manifold $M$ is uniquely determined by 
a bi-Hermitian structure on $M$ whose
torsion form $T_+$ (called {\it flux} 
by physicists) is exact. There is also
a slight generalization of generalized K\"ahler
structures, called {\bf twisted generalized 
K\"ahler structures}, and these are equivalent
to bi-Hermitian structures with arbitrary
(not necessarily exact) torsion form.

In this sense, the notions of generalized K\"ahler 
structure and bi-Hermitian structure are synonymous.

Another notion, also due to physicists, is the notion
of HKT-manifold, which was suggested by Howe and Papadopoulos in 
\cite{HP} and has been much studied since then. 

\begin{definition}
Let $(M,g)$ be a Riemannian manifold, and $I, J, K$ be
complex structures on $(M,g)$ which are Hermitian and
satisfy the quaternionic relations $IJ=-JI=K$.
Then $(M, g, I, J, K)$ is called {\bf a quaternionic
Hermitian hypercomplex manifold}. If, in addition, the Bismut
connections associated to $I$, $J$ and $K$ coincide, then
$(M, g, I, J, K)$ is called {\bf an HKT-manifold}; and if
the Bismut torsion is closed, then $(M, g, I, J, K)$
is called {\bf strong HKT}.
\end{definition}

For more details and examples of hypercomplex manifolds and HKT-geometry, 
please see Section \ref{_HKT_Section_}.

\begin{remark}
An orthogonal connection is uniquely determined by
its torsion (see for example \cite{_Friedrich_Ivanov:physics_}).
Therefore, the Bismut connections associated to
$I$, $J$, $K$ are equal if and only if the
corresponding torsion forms are equal:
\[
I d \omega_I = Jd \omega_J = K d \omega_K.
\]
Consequently, a hypercomplex Hermitian structure $I,J,K$
on a Riemannian manifold $(M,g)$
is HKT with respect to $g$ if and only if the torsion 3-forms 
corresponding to $I$, $J$, $K$ are equal. 
\end{remark}

We finally consider (4,4)-supersymmetry structures 
on Riemannian manifolds. 
These structures were also introduced by Gates, Hull and Ro\v cek in \cite{_GHR_}, 
and can be formulated in Hitchin's and Gualtieri's language as 
{\em generalized hyperk\"ahler structures}. 
These structures were explored in more detail
in \cite{_Bredthauer_}; also see \cite{_Huybrechts_gene_}
and \cite{_Goto_defo_}.

\begin{definition}\label{_4,4_manifold_}
Let $(M,g)$ be a Riemannian manifold, and let $I_+, J_+, K_+$
and $I_-, J_-, K_-$ be two triples of complex structures 
on $(M,g)$ which are hypercomplex Hermitian and HKT with 
respect to $g$.
Denote by $T_+$, $T_-$ the corresponding torsion forms.
Then $(M, g, I_+, J_+, K_+,I_-, J_-, K_-)$ 
is called {\bf a (4,4)-manifold}, or {\bf 
a generalized hyperk\"ahler manifold} if the 
torsion forms $T_\pm$ are closed and satisfy
$T_+=-T_-$.
\end{definition}

\begin{remark}
Note that (4,4)-manifolds are equipped with a plethora
of bi-Hermitian structures. Indeed, 
take a complex structure $V_+$ induced by
the first quaternion action, and a complex structure $U_-$
induced by the second one. Then the
Bismut torsion of $V_+$ is equal to $T_+$, and
the Bismut torsion of $U_+$ is equal to $T_-$, 
implying that $(M,g, V_+, U_-)$ is a bi-Hermitian structure.
\end{remark}

A trivial (and not very interesting) example
of a (4,4)-manifold can be obtained starting
from a hyperk\"ahler manifold $(M,I,J,K,g)$.
Let $I',J',K'$ be another hyperk\"ahler structure
on $(M,g)$; such a triple $I',J',K'$  can be obtained,
for instance, by twisting $I,J,K$ by a
non-zero quaternion. 
Since $M$ is K\"ahler, the Bismut connection
coincides with the Levi-Civita connection, and its torsion 
vanishes. Then $(M, g, I,J,K,I',J',K')$
is a (4,4)-manifold. Clearly, all (4,4)-manifolds
with trivial torsion are obtained this way.

There are not many examples of (4,4)-manifolds
with non-trivial torsion.
In fact, all known examples (except
those provided by Theorem \ref{_4,4_on_insta_Theorem_}) 
are homogeneous or the product of a homogeneous
(4,4)-manifold and a hyperk\"ahler manifold.

Examples of homogeneous (4,4)-manifolds
are not very difficult to obtain. Let $G$
be a semisimple Lie group admitting a 
left-invariant hypercomplex structure
$I_+, J_+, K_+$ (such hypercomplex structures were 
constructed and completely classified
by D. Joyce in \cite{Joyce}). 
Replacing the left multiplication by the right
one, we may also choose a right-invariant
hypercomplex structure $I_-, J_-, K_-$ on 
$G$. The Killing metric $g$ on $G$ is HKT 
(see for example \cite{GrP}), and its
torsion $T_+$ is equal to the fundamental
3-form of $G$. In particular, $T_+$ is closed.
If we replace the left group
multiplication by the right one, then the 
fundamental 3-form of $G$ becomes the opposite of the previous one. 
The torsion form $T_-$ of $I_-, J_-, K_-$ therefore satisfies
$T_-=-T_+$, and $(G,g, I_+, J_+, K_+, I_-, J_-, K_-)$
is a (4,4)-manifold. \\

The main results of the paper are the following.
Consider a bi-Hermitian manifold $X$ of real dimension 4.
It was then shown in \cite{_Hitchin:biherm_} that the
moduli space of anti-self-dual connections on $X$ is also
a bi-Hermitian space. We prove similar results for
hypercomplex, HKT, and (4,4)-structures:

\begin{theorem}\label{_HKT_insta_main_}
Let $(X,I,J,K,g)$ be a compact strong HKT-manifold of real
dimension 4, and $E$ be a smooth complex vector bundle on
$X$. Denote by $\mathcal{M}$ the moduli space 
of gauge-equivalence classes of anti-self-dual 
connections (instantons) on $E$.  Then $\mathcal{M}$ 
is equipped with a natural strong HKT-structure.
\end{theorem}
\begin{proof}
See section \ref{hypercomplex-instanton}.
\end{proof}

\begin{remark}
Compact hypercomplex 4-manifolds were classified in
\cite{Bo}, where it was shown that a compact hypercomplex
4-manifold is either a torus, a K3-surface, 
or a special type of Hopf surface (see section \ref{strong HKT}).
Each of these manifolds admits a strong HKT-structure
(see section \ref{strong HKT}). Therefore, the moduli of 
stable holomorphic $SL(n,\mathbb{C})$-bundles on a given hypercomplex
surface is again hypercomplex.
\end{remark}

{}From Theorem \ref{_HKT_insta_main_}, we deduce the following
theorem.

\begin{theorem}\label{_4,4_on_insta_Theorem_}
Let $(X,I_\pm,J_\pm,K_\pm,g)$ be a compact (4,4)-manifold
of real dimension 4, and $E$ be a smooth complex vector bundle on
$X$. Denote by $\mathcal{M}$ the moduli space 
of gauge-equivalence classes of anti-self-dual 
connections (instantons) on $E$.  Then $\mathcal{M}$ 
is equipped with a natural (4,4)-structure.
\end{theorem}
\begin{proof}
See section \ref{strong HKT}.
\end{proof}


\section{Hypercomplex structures and HKT-metrics}
\label{_HKT_Section_}

\subsection{Hypercomplex manifolds}
\label{hypercomplex manifolds}

A smooth manifold $M$ equipped with three complex structure operators 
$I, J, K: TM \rightarrow TM$
that satisfy the quaternionic identities 
\begin{equation}\label{quaternionic identities} 
IJ = -JI = K 
\end{equation}
is said to be {\em hypercomplex} or to admit a {\em hypercomplex structure}.
The complex structures $I$, $J$, and $K$ induce other almost 
complex structures on $M$ 
of the form $L := aI + bJ + cK$ 
for all real numbers $a,b,c$ such that $a^2+b^2+c^2=1$; 
that these almost complex structures are in fact integrable follows from Obata
\cite{Obata, Kaledin}.
Given such a complex structure $L$ on $M$,
we will denote by $(M,L)$ the manifold $M$ considered as a 
complex manifold with respect to $L$.

In this paper, we study the moduli spaces of instantons (solutions to the anti-self-dual
Yang-Mills equations) on compact hypercomplex 
4-manifolds; we show, in particular, that these moduli spaces 
admit a natural hypercomplex structure which is induced from the hypercomplex 
structure on the 4-manifolds.

Compact hypercomplex 4-manifolds were classified by Boyer
who showed that if $(X,I,J,K)$ is a compact hypercomplex 4-manifold, 
then $X$ is either a torus, a $K3$ surface, or a quaternionic Hopf surface
(see \cite{Bo}, Theorem 1).
Recall that a {\em quaternionic Hopf surface} $X$ can be defined as the quotient 
of the non-zero quaternions $\mathbb{H} - \{0\}$ 
by a cyclic group generated by some $q \in \mathbb{H}$ with $|q| > 1$, 
where $\langle q \rangle$ acts on $\mathbb{H} - \{0\}$ by right multiplication:
\begin{equation}\label{Hopf surface} 
X := (\mathbb{H} - \{0\}) / \langle q \rangle.
\end{equation} 
Note that the action of left multiplication by 
$i$, $j$, and $k$ commutes with the action of right multiplication by $q$;
hence, the hypercomplex structure $\{ I,J,K \}$ on $\mathbb{H}$ induced 
by left multiplication by $i,j,k$, respectively, 
descends to a hypercomplex structure on the Hopf surface. 
Furthermore, any quaternionic Hopf surface 
in Ma. Kato's classification (\cite{K}, Proposition 8)
is isomorphic to a finite cover of a Hopf surface of the form \eqref{Hopf surface},
thus acquiring the same hypercomplex structure from $\mathbb{H}$.

It has been known for some time that instanton moduli spaces on tori and K3 surfaces
admit hypercomplex structures.
In this article, we show that
this is true for hypercomplex 
Hopf surfaces; this is done by identifying the instanton moduli spaces
with moduli spaces of stable bundles, implying, in particular, that we will 
consider metrics on these 4-manifolds which are Hermitian\footnote
{A Riemannian metric $g$ on a smooth manifold $M$ with complex structure $L$ is called 
{\em Hermitian} if $g(LX, LY) = g(X,Y)$ for all vector fields $X$ and $Y$ on $M$.} 
with respect to every complex structure (see section \ref{stable bundles}).

In \cite{Joyce}, page 747, D. Joyce suggested that 
the space of instantons on quaternionic Hopf surfaces
can be obtained through quaternionic reduction. 
Similar results were obtained independently by 
Oliver Nash  and Gil Cavalcanti in unpublished papers
\cite{_Oliver_Nash_} and 
\cite{_Cavalcanti_Redu_}, using the methods of hypercomplex reduction
(Nash) and reduction of Courant algebroids applied to
generalized K\"ahler geometry (Cavalcanti).

A Riemannian metric $g$
on a hypercomplex manifold $(M,I,J,K)$
is called {\em hyperhermitian} or {\em quaternionic Hermitian} if it is 
Hermitian
with respect to every complex structure $L$ on $M$ induced by $I,J,K$.
In addition, if a hyperhermitian metric $g$ is K\"ahler\footnote
{One can associate to any Hermitian metric $g$ on $(M,L)$ the 2-form 
$\omega_L(\cdot,\cdot) := g(L\cdot,\cdot)$, called the {\em Hermitian form} of $g$.
A Hermitian metric $g$ is then said to be {\em K\"ahler} if its Hermitian form $\omega_L$ is $d$-closed.} 
for all complex structures on $M$, then it is called {\em hyperk\"ahler};
the Euclidean metric on $\mathbb{H}$ is an example of a hyperk\"ahler metric.
Note that hyperhermitian metrics exist on all hypercomplex manifolds $(M,I,J,K)$; 
indeed, one can construct a hyperhermitian metric on $M$ by taking any 
Riemannian metric on $M$ and averaging it over the natural $SU(2)$-action on $M$ 
(induced by multiplication by the quaternions). 
However, hyperk\"ahler metrics only exist if the underlying manifold admits 
K\"ahler metrics; for instance, quaternionic Hopf surfaces do not admit 
K\"ahler metrics since they have odd first Betti number, implying that they do not 
admit K\"ahler structures.

One can endow tori and $K3$ surfaces with hyperk\"ahler metrics 
(for details, see \cite{_Besse:Einst_Manifo_});
quaternionic Hopf surfaces are therefore the only compact hypercomplex 
4-manifolds on which hyperhermitian metrics are never hyperk\"ahler.
One can nonetheless construct hyperhermitian metrics on quaternionic Hopf surfaces 
which are Gauduchon\footnote
{A Hermitian metric $g$ on an $n$-dimensional complex manifold $(M,L)$ 
is called Gauduchon if the $(n-1)$-th power
its Hermitian form $\omega_L$ is $dd^c_{L}$-closed,  
where $d^c_{L}$ is the twisted differential which acts as 
$(-1)^m L \circ d \circ L$ on $m$-forms. 
Note that although K\"ahler metrics are Gauduchon, the converse is in general not true.}  
with respect to every complex structure. Consider, for instance, 
a quaternionic Hopf surface 
of type \eqref{Hopf surface}.  
Let $r$ be the Euclidean length on $\mathbb{H}$ and let $\varphi := r^2$. 
The 2-forms 
\[ \omega_L := \frac{dd^c_L\varphi}{\varphi}, \]
where $d^c_L$ denotes the twisted differential,
are then $\langle q \rangle$-invariant, and descend to 2-forms on $X$ which induce the same metric $g$ on $X$,
that is, $g(\cdot,\cdot) = \omega_L(\cdot,L\cdot)$ for all complex structures $L$.
The metric $g$ is thus hyperhermitian.
Moreover, a direct computation shows that 
\[ d^c_I \omega_I = d^c_J \omega_J = d^c_K\omega_K = H,\]
where $H$ is a $d$-closed 3-form, implying that
$g$ is Gauduchon with respect to every complex structure on $X$
induced by $I,J,K$. 

\subsection{HKT-metrics and (4,4)-symmetry}
\label{strong HKT}

Consider a hypercomplex manifold $(M,I,J,K)$. A hyperhermitian metric $g$ on $M$
is  then called an {\em HKT-metric} if
\[d^c_I \omega_I = d^c_J \omega_J = d^c_K\omega_K = H,\]
for some 3-form $H$,
where $\omega_L$ is the Hermitian form of $g$ and $d^c_L$ is the twisted
differential, corresponding to the complex structures $L= I, J, K$.
Moreover, if $H$ is $d$-closed, then $g$ is said to be a {\em strong HKT-metric}.
Note that for any complex structure $L$ on $M$, 
the skew-symmetric torsion of the Bismut connection on $(M,L)$
is equal to the 3-form $-2H$ (see \cite{HP}, or \cite{GrP} Proposition 1).
Furthermore, an HKT-metric $g$ is hyperk\"ahler if and only if $H=0$
(hyperk\"ahler metrics are in fact strong HKT). However,
on a manifold that does not admit K\"ahler metrics, one has $H \neq 0$, 
hence the terminology HKT which stands for 
``hyperk\"ahler metric with torsion''.

There exists another characterisation of HKT-metrics.
Let $g$ be a hyperhermitian metric on $(M,I,J,K)$ and 
let $\omega_I$, $\omega_J$, and $\omega_K$ be its 
Hermitian forms for $I$, $J$, and $K$, respectively.
Set 
\[ \Omega := \omega_J + \sqrt{-1}\omega_K. \]
Then $\Omega$ is a $(2,0)$-form on $(M,I)$, which can be used to determine
whether the metric is HKT. 
Indeed, one can show 
that the metric $g$ is HKT if and only if 
it satisfies the condition $\partial \Omega = 0$,
where $\partial = \frac{1}{2}(d + \sqrt{-1}d^c_I)$
(see \cite{HP} and \cite{GrP}, Proposition 2).
Consequently, since $\partial \Omega$ is a $(3,0)$-form on $(M,I)$, if $M$ is a 4-manifold,
any hyperhermitian structure is an HKT-structure.
This implies that on hypercomplex 4-manifolds, strong HKT-metrics  
are equivalent to hyperhermitian metrics that are Gauduchon 
with respect to all complex structures. 
Every hypercomplex compact 4-manifold therefore admits a strong HKT-metric:
referring to section \ref{hypercomplex manifolds}, 
tori and K3 surfaces admit hyperk\"ahler metrics, 
and quaternionic Hopf surfaces admit metrics which are
Gauduchon with respect to every complex structure.

Let us now consider the quaternionic Hopf surface 
\[ X:= (\mathbb{H} - \{ 0 \}) / \langle q \rangle\]
with $q \in \mathbb{R}$.
One can then endow $X$ with two natural hypercomplex structures.
We have seen that left multiplication by $i$, $j$, and $k$ on $\mathbb{H}$
induces a hypercomplex on $X$, which we now denote $I_+,J_+,K_+$. 
The other hypercomplex structure on $X$ corresponds to right multiplication
by $i$, $j$, and $k$ on $\mathbb{H}$
(since $q$ is real, its action on $\mathbb{H}$ commutes with the action of
right multiplication by $i$, $k$, and $k$); we will denote this second hypercomplex
structure $I_-,J_-,K_-$.
Note that any hypercomplex structure on $X$ induced by one of these two hypercomplex 
structures is orientation preserving. 
Moreover, one can verify that the 2-forms
$(dd^c_L \varphi)/ \varphi$, where $\varphi = r^2$ and $r$ is 
the Euclidean length on $\mathbb{H}$, induce the same metric on $X$ which
is a strong HKT-metric for both hypercomplex structures 
$I_+,J_+,K_+$ and $I_-,J_-,K_-$.
In fact,
\[d^c_{I_+} \omega_{I_+} = d^c_{J_+} \omega_{J_+} = d^c_{K_+}\omega_{K_+} = H,\]
and
\[d^c_{I_-} \omega_{I_-} = d^c_{J_-} \omega_{J_-} = d^c_{K_-}\omega_{K_-} = -H,\]
for some $d$-closed 3-form $H$.
Finally, the two families of hypercomplex structures are independent, 
in the sense that no complex structure
induced by $I_+,J_+,K_+$ can be written as a linear combination of $I^-,J^-,K^-$,
and vice-versa. Hence, every pair $(L_+,L_-)$ with $L=I$, $J$, or $K$ defines a bi-Hermitian
structure $X$. 
A pair of strong HKT structures that satisfy the above properties is known as a
$(4,4)$-structure (see Definition \ref{_4,4_manifold_}). 
The Hopf surface endowed with its two hypercomplex structures
$I_+,J_+,K_+$ and $I_-,J_-,K_-$ is then an example of $(4,4)$-symmetry.

It then follows from Theorems \ref{_HKT_insta_main_} and \ref{_HKT_moduli_} that
the natural $L^2$-metric $g_{L^2}$ on the instanton moduli space $\mathcal{M}$ 
is strong HKT with respect to the hypercomplex structures induced by both 
$I_+,J_+,K_+$ and $I_-,J_-,K_-$. 
Moreover, a result of Hitchin's \cite{_Hitchin:biherm_} 
shows that each pair $(L_+,L_-)$ induces a bi-Hermitian structure
on $\mathcal{M}$ with respect to the $L^2$ metric $g_{L^2}$.  This
implies that the instanton moduli space also 
admits a $(4,4)$-structure.


\section{Instanton moduli spaces}
\label{_Instanton_Section_}


\subsection{Hermitian-Einstein connections and stable bundles}
\label{stable bundles}

Let $X$ be a compact complex surface with fixed Gauduchon metric $g$ and 
Hermitian form $\omega$; in particular, we have $\partial\bar{\partial}\omega = 0$.
Consider a smooth vector bundle $E$ on $X$ and let $h$ be a Hermitian metric in $E$. 
The space of $h$-unitary connections in $E$ is denoted $\mathcal{A}(E,h)$.

Recall that a connection $A$ in $\mathcal{A}(E,h)$ is called 
{\em $g$-Hermitian-Einstein} 
if its curvature 2-form $F_A$ is of type $(1,1)$ and satisfies
\[ \sqrt{-1}\Lambda_g F_A = \gamma_A \cdot \rm{id}_E \]
for some $\gamma_A \in \mathbb{R}$, where $\Lambda_g$ is the contraction of 2-forms by $\omega$.
Note that all Hermitian-Einstein connections are integrable and therefore induce 
holomorphic structures in $E$. Moreover, irreducible $g$-Hermitian-Einstein connections
give rise to $g$-stable holomorphic structures in $E$, where $g$-stability is defined as follows. 

Stability with respect to Gauduchon metrics is an 
extension of Mumford-Takemoto stability and thus requires the notion of degree.
Given the Gauduchon metric $g$ on $X$, the {\it degree} of a holomorphic line bundle $L$ on $X$ 
is defined, up to a multiplicative constant, by
\[ \deg{L} := \frac{1}{2\pi}\int_M F \wedge \omega,\]
where $F$ is the curvature of a Hermitian connection on $L$, compatible with 
$\bar{\partial}_L$. Since any two such forms $F$ differ by a $\partial \bar{\partial}$-exact form
and $\partial \bar{\partial} \omega = 0$, the degree does not depend on the choice of connection.

Note that flat line bundle have degree zero since the curvature of any connection 
on such bundles is zero; in particular, the trivial line bundle has degree zero. 
Furthermore, if the metric $g$ is K\"ahler, then we get the usual topological degree; 
otherwise, the degree is 
not a topological invariant, as it can take continua of values in $\mathbb{R}$.

The {\it degree} of a torsion-free coherent sheaf $\mathcal{E}$ on $X$ is given by
\[ \deg(\mathcal{E}) := \deg(\det{\mathcal{E}}), \]
where $\det{\mathcal{E}}$ is the determinant line bundle of $\mathcal{E}$,
and the {\it slope of $\mathcal{E}$} is defined as
\[ \mu(\mathcal{E}) := \deg(\mathcal{E})/\text{rk}(\mathcal{E}).\]
A torsion-free coherent sheaf $\mathcal{E}$ on $X$ is then said to be {\em $g$-(semi)stable} 
if and only if for every proper coherent subsheaf $\mathcal{S} \subset \mathcal{E}$  
we have 
\[ \mu(\mathcal{S}) \leq \mu(\mathcal{E}), \]
with strict inequality for $g$-stable bundles.
One can then show that a holomorphic vector bundle $\mathcal{E}$ is $g$-stable if and only if it admits
an irreducible $g$-Hermitian-Einstein connection; this was done by Buchdahl \cite{Buchdahl} for surfaces
and Li and Yau \cite{LY} for all Hermitian manifolds. 
The one-to-one correspondence between irreducible Hermitian-Einstein 
connections and stable holomorphic structures in a smooth vector bundle
is known both as the Donaldson-Uhlenbeck-Yau correspondence 
and the Kobayashi-Hitchin correspondence; 
a comprehensive reference on the subject is the book \cite{LT}.
\vspace{.05in}

Recall that a
connection $A$ on $E$ is called {\em anti-self-dual} (ASD) if its curvature
$F_A$ satisfies the equation: 
\[ \ast F_A = -F_A, \] 
or equivalently, if $F_A$ is a matrix of anti-self-dual 2-forms.
Since the anti-self-dual forms on a Hermitian 4-manifold $(M,J,g)$ 
consist of $(1,1)$-forms which are orthogonal to the Hermitian form $\omega$ of $g$,
a connection $A$ is ASD if and only if it is $g$-Hermitian-Einstein
and $\Lambda_g F_A = 0$.
Irreducible anti-self-dual connections in $E$ therefore induce $g$-stable holomorphic structures 
of degree zero in $E$. 

Let us now consider a compact hypercomplex 4-manifold $(X,I,J,K)$ equipped with a strong 
HKT-metric $g$. Let $E$ be a smooth complex vector bundle on $X$
and $h$ be a Hermitian metric in $E$. An $h$-unitary connection $A$ in $E$
is said to be 
{\em hyperholomorphic} if it is integrable with respect to every complex
structure $L$ on $X$.
Note that anti-self-dual forms on $X$ are of type $(1,1)$ 
with respect every complex structure $L$ on $X$ induced by $I,J,K$. 
Anti-self-dual connections in $E$ are therefore hyperholomorphic.
The converse is also true:

\begin{theorem}
Let $E$ be a vector bundle with Hermitian metric $h$ 
on a compact strong HKT 4-manifold $(X,I,J,K,g)$.
Then, an $h$-unitary connection $A$ in $E$
is hyperholomorphic if and only if $A$ is anti-self-dual.              
\end{theorem}
\begin{proof}
See for example \cite{V-hyper}, sections 1 and 2, for a proof in the case where
$(X,I,J,K)$ admits a hyperk\"ahler metric.
The arguments used in \cite{V-hyper} extend, however, to all hypercomplex 4-manifolds. 
\end{proof}

Consequently, since anti-self-dual forms on a hyperhermitian 4-manifold $(X,I,J,K,g)$ are
orthogonal to the Hermitian form $\omega_L$ of $g$ for all complex structures $L$ on $X$,
we see that a connection in $E$ induces a $g$-stable holomorphic structure in $E$
with respect to all complex structures on $X$
if and only if it is anti-self-dual. 

In the next section, we study moduli spaces of connections on compact strong HKT 
4-manifolds $(X,I,J,K,g)$.
We show in particular that these moduli spaces admit hypercomplex structures; 
this is done by identifying these moduli spaces with moduli spaces of $g$-stable holomorphic 
structures, for all complex structures on the 4-manifolds. We therefore only consider
connections that induce $g$-stable 
holomorphic structures for all complex structures on $X$, i.e., 
anti-self-dual connections.

\subsection{Hypercomplex structures and HKT-metrics}
\label{hypercomplex-instanton}

Let $(X,I,J,K,g)$ be a compact strong HKT 4-manifold
and let $E$ be a smooth complex vector bundle on $X$.
The space $\mathcal{A}^{ASD}$ of all irreducible anti-self-dual connections (instantons) in $E$
is then a $\mathcal{G}$-principal bundle,
where $\mathcal{G}$ is the group of gauge transformations of $E$.
The quotient space 
\[ \mathcal{M} := \mathcal{A}^{ASD}/\mathcal{G} \]
is the moduli space of gauge equivalence classes of anti-self-dual
connections in $E$.
Moreover, the $L^2$ metric
\begin{equation}\label{metric}
(a_1,a_2) = - \int_M {\tt tr}(a_1 \wedge \ast a_2)
\end{equation}
on $\mathcal{A}^{ASD}$ descends to a metric $g_{L^2}$ on the moduli space $\mathcal{M}$.

Referring to the previous section,
instantons correspond to connections in $E$ that are integrable 
with respect to every complex structure on $X$ induced by the quaternions.
The moduli space $\mathcal{M}$ can therefore be identified
via the Kobayashi-Hitchin correspondence with the moduli space
$\mathcal{M}^{st}_L$ of isomorphism classes of $g$-stable holomorphic 
structures in $E$, for any complex structure $L$ on $X$.
The moduli space $\mathcal{M}$ thus inherits a natural complex structure from 
$\mathcal{M}^{st}_L$, for every complex structure $L$ on $X$, 
which can be described as follows.

Recall that the tangent space to the moduli space $\mathcal{M}$ at any point $[A]$ 
can be identified with the horizontal 
subspace at $A$ of any connection on the principal $\mathcal{G}$-bundle  
\[\mathcal{P}:= \mathcal{A}^{ASD} \rightarrow \mathcal{M}.\]
Moreover, since the difference between any two connections in $E$ is a 1-form with values 
in $sl(E)$, where $sl$ denotes trace-free endomorphisms, 
then every element of $\mathcal{A}^{ASD}$ is of the form $A + a$ for some $a \in A^1(sl(E))$. 
Suppose that one fixes a complex structure $L$ on $X$.
The horizontal subspace at $A$ is then chosen to be the set of 1-forms 
$a$ such that $\Lambda_g d^c_L a = 0$
(where the subscript $A$ in $d^c_L$ is suppressed for clarity),
whereas the vertical subspace is the tangent space of the $\mathcal{G}$-orbit through $A$,
giving us the following local model:
\begin{equation}\label{tangent space - instantons} 
T_{[A]}\mathcal{M} = 
\{ a \in A^1(sl(E)) \ | \ \mbox{$d_A^+ a = 0$ and $\Lambda_gd_L^c a = 0$} \}
\end{equation}

\noindent
(see for example \cite{LT} for more details).

The advantage of using this particular connection on the 
$\mathcal{G}$-bundle $\mathcal{P}$ is twofold.
The complex structure $\tilde{L}$ on $\mathcal{M}$
induced from the natural complex structure on $\mathcal{M}^{st}_L$
has a very simple expression at any given point $[A]$.
Indeed, note that the complex structure $L$ on $X$ decomposes $A^d(sl(E))$ into
components $A^{p,q}(sl(E))$ with $p+q=d$;
given this decomposition, the complex structure $\tilde{L}$ 
on $\mathcal{M}$ is the operator
\begin{equation}\label{complex structure} 
\tilde{L}(a) := \sqrt{-1}\left(a^{0,1} - a^{1,0}\right),
\end{equation} 
for any $a \in T_{[A]}\mathcal{M}$.
Furthermore, the metric $g_{L^2}$ on $\mathcal{M}$ is Hermitian with respect to $\tilde{L}$,
and has the following properties:
\begin{theorem}[L\"{u}bke-Teleman]\label{Lubke-Teleman}
The natural $L^2$ metric $g_{L^2}$ on $(\mathcal{M},\tilde{L})$ is Hermitian 
and its Hermitian form $\tilde{\omega}_{\tilde{L}}$ is such that:

(i) Let $\theta$ be the curvature of the connection on the principal $\mathcal{G}$-bundle
$\mathcal{P}$ which has horizontal subspaces \eqref{tangent space - instantons}, 
and denote by $\bar{a}_i$ any horizontal lift of 
$a_i \in T_{[A]}\mathcal{M}$ to $\mathcal{A}^{ASD}$. 
Then,
\[ \tilde{\omega}_{\tilde{L}}(a_1,a_2) = \int_M \omega_L \wedge {\tt tr}(\bar{a}_1 \wedge \bar{a}_2), \]
\noindent
and
\begin{equation}\label{differential}
d^c_{\tilde{L}} \tilde{\omega}(a_1,a_2,a_3)  = \frac{1}{3} \sum_{\sigma \in S_3} 
(-1)^{{\tt sign} \sigma} \int_M d^c_L \omega \wedge 
{\tt tr}\left(\theta(\bar{a}_{\sigma(1)},\bar{a}_{\sigma(2)})\bar{a}_{\sigma(3)}\right).
\end{equation}

(ii) $dd^c_{\tilde{L}} \tilde{\omega}_{\tilde{L}} = 0$.
\end{theorem} 
\begin{proof}
For details see \cite{LT}, Theorem 5.3.6 and Lemma 5.3.7.
\end{proof}

\begin{remarks}\label{independence}
Given any element $a \in A^1(sl(E))$ we have:
\[ d^\ast_A a = \Lambda_g d^c_L a + \ast (d^c_L\omega_L \wedge a). \]
The horizontal slices \eqref{tangent space - instantons} can then be described as
\begin{equation}\label{horizontal slice}
\{ a \in A^1(sl(E)) \ | \ \mbox{$d_A^+ a = 0$ and $d^\ast_A a = \ast (d^c_L\omega_L \wedge a)$} \},
\end{equation}
\noindent
which implies the following.
 
(i) If the metric $g$ is K\"ahler with respect to $L$, then we have
\[ d^\ast_A a = \Lambda_g d^c_L a \]
since $d^c_L \omega = 0$;
in this case \eqref{horizontal slice} reduces to 
\[\{ a \in A^1(sl(E)) \ | \ \mbox{$d_A^+ a = 0$ and $d^\ast_A a = 0$} \},\]
so that
\eqref{tangent space - instantons} is the usual local model for the tangent space to $\mathcal{M}$
(namely the orthogonal complement in $\mathcal{A}^{ASD}$ to the tangent space 
of the gauge orbit at $A$, with respect to the $L^2$ metric \eqref{metric}).

(ii) If the metric $g$ is strong HKT, then 
\begin{equation}\label{differentials}
d^c_I \omega_I = d^c_J \omega_J = d^c_K \omega_K = H
\end{equation}
for some $d$-closed 3-form $H$ on $X$.
Consequently, given the description \eqref{horizontal slice} of the 
horizontal spaces, we see that the tangent space to $\mathcal{M}$ at $[A]$
is the same for all complex structures $L$; 
one can therefore compose the complex structures $\tilde{L}$ on $\mathcal{M}$.
Moreover, our choice of connection on $\mathcal{P}$ is independent of the complex structure,
so that its connection matrix $\theta$ is the same for all complex structures
$L$. This, combined with \eqref{differentials}, \eqref{differential}, and 
Theorem \ref{Lubke-Teleman} (ii), gives us that
\[ d^c_{\tilde{I}} \omega_{\tilde{I}} = d^c_{\tilde{J}} \omega_{\tilde{J}} = 
d^c_{\tilde{K}} \omega_{\tilde{K}}=\tilde{H}, \]
where $\tilde{H}$ is a $d$-closed 3-form.
\end{remarks} 

Referring to Remark \ref{independence} (i), the complex structures 
$\tilde{I}$, $\tilde{J}$, and $\tilde{K}$ on $\mathcal{M}$
can be composed; moreover, these complex structures satisfy the quaternionic
identities \eqref{quaternionic identities}, since we have the following:

\begin{lemma}
The complex structures $I$ and $J$ on $X$ induce
complex structures $\tilde{I}$ and $\tilde{J}$
on $\mathcal{M}$ that anti-commute.
\end{lemma}
\begin{proof}
A section $a$ of $A^1(sl(E))$ can be written locally as
\[ a = \Sigma a_i \otimes s_i, \]
with $a_i \in A^1(X)$ and $s_i \in sl(E)$. For any complex structure $L$
on $X$, one therefore has
\[ \tilde{L}(a) = - \Sigma L(a_i) \otimes s_i, \]
which is independent of the local trivialisation.
Hence, since $IJ=-JI$, we have that $\tilde{I}\tilde{J}(a) = - \tilde{J}\tilde{I}(a)$
for all $a \in A^1(sl(E))$.
\end{proof}

We therefore have the following results:
\begin{theorem}
Let $(X,I,J,K)$ be a compact hypercomplex 4-manifold  and let $g$ be a 
strong HKT-metric on $X$.
Let $E$ be a fixed smooth complex vector bundle on $X$.
The moduli space $\mathcal{M}$
of gauge-equivalence classes of anti-self-dual connections on $E$
then admits a hypercomplex structure. \hfill $\square$
\end{theorem}
 
\begin{theorem}\label{_HKT_moduli_}
Let $(X,I,J,K)$ be a compact hypercomplex 4-manifold and let $g$ 
be a strong HKT-metric on $X$.
Fix a smooth complex vector bundle $E$ on $X$ and 
consider the moduli space $\mathcal{M}$
of gauge-equivalence classes of anti-self-dual connections on $E$. 
The $L^2$ metric on $\mathcal{M}$ is then strong HKT. 
\end{theorem}
\begin{proof}
This follows from Theorem \ref{Lubke-Teleman} and Remark \ref{independence} (ii).
\end{proof}

{\bf Acknowledgements:}
We are grateful to M. Gualtieri for insightful advice
and comments.

\hfill

\noindent{\sc R. Moraru\\
Department of Pure Mathematics, University of Waterloo\\
200 University Avenue West,Waterloo, ON, Canada N2L 3G1.}\\\\
\ \\
\noindent {\sc Misha Verbitsky\\
University of Glasgow, Department of Mathematics, \\
15 University Gardens, Glasgow G12 8QW, Scotland.}\\
\   \\
{\sc  Institute of Theoretical and
Experimental Physics \\
B. Cheremushkinskaya, 25, Moscow, 117259, Russia }\\
\  \\
\tt verbit@maths.gla.ac.uk, \ \  verbit@mccme.ru

\end{document}